\documentclass[reprint,superscriptaddress,amsmath,amssymb,aps,nofootinbib,prx]{revtex4-1}
\usepackage{amsthm}
\usepackage[colorlinks=true,urlcolor=blue,citecolor=blue,linkcolor=blue]{hyperref}
\usepackage[shortlabels]{enumitem}
\usepackage{amsfonts}
\usepackage{graphicx}
\usepackage[caption = false]{subfig}
\usepackage{multirow}
\usepackage{float}
\usepackage{bbold}
\usepackage{stmaryrd}
\usepackage{color}
\usepackage{ragged2e}
\usepackage{hyperref}
\usepackage{verbatim}
\usepackage{cases}
\usepackage{amsfonts}
\usepackage{amssymb}
\usepackage[]{mathrsfs}
\usepackage[table]{xcolor}
\usepackage{epsfig}
\usepackage{bm}
\usepackage{setspace}
\usepackage{enumerate}
\usepackage{dcolumn}% Align table columns on decimal point
\usepackage{bm}% bold math
\usepackage{enumitem}
\usepackage{caption} 
\usepackage{tikz}
\usepackage{multirow, array} % para las tablas

\usepackage{graphicx}% Include figure files
\usepackage{dcolumn}% Align table columns on decimal point
\usepackage{bm}% bold math
\usepackage{color}
\usepackage{amssymb}
\usepackage{lipsum}
\usepackage{amsmath}
\usepackage{graphicx}
\usepackage{array}
\usepackage{booktabs}
\usepackage[english]{babel}
\begin{document}

%\preprint{APS/123-QED}

\title{ Data-driven reconstruction of chaotic dynamical equations: the Hénon-Heiles type system}% Force line breaks with \\
%\footnote{explicar en nota a pie de pagina USMER y sector}
%\author{Author}

\author{A. M. Escobar-Ruiz}%
\affiliation{Departamento de F\'{i}sica, Universidad Aut\'onoma Metropolitana Unidad Iztapalapa, San Rafael Atlixco 186, 09340 Cd. Mx., M\'exico}

\author{L. Jiménez-Lara}%
\affiliation{Departamento de F\'{i}sica, Universidad Aut\'onoma Metropolitana Unidad Iztapalapa, San Rafael Atlixco 186, 09340 Cd. Mx., M\'exico}

\author{P. M. Juárez-Florez}%
\affiliation{Departamento de F\'{i}sica, Universidad Aut\'onoma Metropolitana Unidad Iztapalapa, San Rafael Atlixco 186, 09340 Cd. Mx., M\'exico}

\author{F. Montoya-Molina}%
\affiliation{Departamento de F\'{i}sica, Universidad Aut\'onoma Metropolitana Unidad Iztapalapa, San Rafael Atlixco 186, 09340 Cd. Mx., M\'exico}

\author{J. Moreno-Sáenz}%
\affiliation{ALITECS Corporation, Tsurumi 74-1, 2300046, Yokohama, Japan,}

\author{M. A. Quiroz-Juarez}%
\email{maqj@fata.unam.mx}
\affiliation{Centro de F\'{i}sica Aplicada y Tecnolog\'{i}a Avanzada, Universidad Nacional Aut\'onoma de M\'exico, Boulevard Juriquilla 3001, 76230 Quer\'{e}taro, M\'exico}


%\date{August 22, 2022}% It is always \today, today,

\begin{abstract}
In this study, the classical two-dimensional potential $V_N=\frac{1}{2}\,m\,\omega^2\,r^2 +   \frac{1}{N}\,r^N\,\sin(N\,\theta)$, $N \in {\mathbb Z}^+$, is considered. At $N=1,2$, the system is superintegrable and integrable, respectively, whereas for $N>2$ it exhibits a richer chaotic dynamics. For instance, at $N=3$ it coincides with the Hénon-Heiles system. The periodic, quasi-periodic and chaotic motions are systematically characterized employing time
series, Poincaré sections, symmetry lines and the largest Lyapunov exponent as a function of the energy $E$ and the
parameter $N$. Concrete results for the lowest cases $N=3,4$ are presented in complete detail. This model is used as a benchmark system to estimate the accuracy of the Sparse Identification of Nonlinear Dynamical Systems (SINDy) method, a data-driven algorithm which reconstructs the underlying governing dynamical equations. We pay special attention at the transition from regular motion to chaos and how this influences the precision of the algorithm. In particular, it is shown that SINDy is a robust and stable tool possessing the ability to generate non-trivial approximate analytical expressions for periodic trajectories as well. 
\end{abstract}

\keywords{}%Use showkeys class option if keyword
\maketitle

%\tableofcontents

%\section{\label{sec:level1}First-level heading:\protect\\ The line
%break was forced \lowercase{via} \textbackslash\textbackslash}
%\section{\label{sec:level1}Introducción}

\section{Introduction} 

The interest on the Hénon-Heiles potential \cite{Henon} was originally motivated by the astronomy community to investigate the 3D motion of a moving star in the gravitational field of an axis-symmetric galaxy. From the two evident integrals, namely the energy and the angular momentum around the symmetry axis, astronomers pondered the question on the existence of an additional conserved quantity $\cal I$ in the Liouville sense (the so called \textit{third integral}) that would make the system integrable. This searching, without any knowledge of the physical origin of $\cal I$ or about its possible mathematical expression. It turned out that the system is chaotic and no third integral exists. Accordingly, several analytical and numerical systematic studies on galactic potentials were conducted in the 1960s (see, for instance Contopoulos \cite{contopoulos1966}) whilst recent references on the closely related classical and quantum chaos can be found in  \cite{gutzwiller2013chaos},\cite{contopoulos2002order}.

Exploiting the conservation of angular momentum, Hénon and Heiles eventually obtained a reduced  Hamiltonian with two degrees of freedom only. This reduced problem is called the Hénon-Heiles system (HHS). The corresponding classical Hamiltonian is of the form:
\begin{equation}
\label{HenHei}
{H} \ = \  \frac{1}{2}\big( p_x^2 \ + \ p_y^2 \big) \ + \  \frac{1}{2}\left(x^2+y^2\right) \ + \ g\,\big(\,y\,x^2 \ -\ \frac{1}{3}\,y^3 \,\big)\ , 
\end{equation}
where $x$ is the altitude and $y$ denotes the radius (the distance to the aforementioned symmetry axis), $p_x$ and $p_y$ are the corresponding canonical momenta and $g>0$, the coupling constant, is a real parameter. The total energy $H=E$ of the above 2D system is the single conserved quantity of the problem. The potential of HHS can be viewed as a two-dimensional isotropic harmonic oscillator perturbed with nonlinear cubic terms. 
Also, it is worth mentioning that the Hamiltonian (\ref{HenHei}) possesses a dihedral symmetry, i.e., a $D_3$ symmetry. As indicated in \cite{Henon}, the Hamiltonian (\ref{HenHei}) is one of the simplest models to capture the relevant dynamics of non-integrable systems. In particular, it incorporates the physical feature that stars are trapped by an attractive potential, and those stars with sufficient energy may escape from the galaxy.

A straightforward numerical analysis \cite{Henon} of HHS revealed the existence of a critical value $E_c$ of energy separating the bounded trajectories from the unbounded ones. For small values $E \ll E_c$, the cubic terms in the potential are negligible and the system displays a regular behaviour mostly. As the energy increases, these terms become relevant and a mixed combination of regular and chaotic motion occurs. Eventually, at $E \sim E_c$ the chaotic regions dominate the phase space landscape.  

Various generalizations of the HHS have been studied in the literature either by considering the coefficients of the cubic terms, in the potential, as free parameters or by adding higher order terms. There is a plethora of works about the integrability or non-integrability of Hénon-Heiles type systems with two parameters; see, for example, Llibre and  Jiménez \cite{llibre2011periodic} and references included therein. Remarkably, the integrability of a generalized HHS that introduces two parameters in the quadratic terms of the potential was demonstrated by Grammaticos et al. \cite{GRAMMATICOS1982111}, see also \cite{FORDY1991204, Conte2005}.

Needless to say that several applications of Hénon-Heiles type systems appear in different fields of physics. For example, an interesting theoretical study of a generalized 3D Hénon-Heiles potential was applied to model the dynamics of ions in a 3D axialy symmetric Penning trap \cite{lanchares2002perturbed},\cite{horvath1998ion}.

In the present study we analyze a one-parametric generalization of (\ref{HenHei}) with the following properties: (I) the potential is the sum of an isotropic 2D harmonic oscillator plus a polynomial function of degree $N=1,2,\,\ldots,$ in the variables $(x,y)$, (II) the system admits a 2D discrete $D_N$ symmetry (group of symmetries of a regular polygon with $N$ sides), (III) there exists a critical value of energy that separates the bounded and unbounded motion, and (IV) mixed regions in the phase space of regular and chaotic dynamics appears for $N\geq 3$. 

The goal of this work is two-fold. Firstly, we aim to investigate the rich dynamics of the system as a function of $N$ and the energy $E$. To this end, we will employ time series, Poincaré sections, symmetry lines and the largest characteristic Lyapunov exponents. For instance, in the computation of the symmetry lines we exploit the discrete rotational symmetry of the Hamiltonian. 
Secondly, we will take this model as a benchmark system to test qualitative and quantitatively the algorithm of Machine Learning SINDy which \textit{reconstructs} the governing dynamical equations of a given system using only the time series data. In particular, we are interested in the influence of the degree of chaoticity on the accuracy provided by SINDy.

\section{The model}

We consider a classical non-relativistic system with two degrees of freedom on the plane. The Hamiltonian is given by
\begin{equation}
\label{HN}
{\cal H}_{{}_N}\ = \  \frac{1}{2\,m}\big( p_r^2 \, + \, \frac{1}{r^2}\,p_{\theta}^2 \big) \  + \ V_N(r,\,\theta)  \ , 
\end{equation}
with scalar potential
\begin{equation}
\label{VN}
V_N(r,\,\theta) \ = \ \frac{1}{2}\,m\,\omega^2\,r^2 \ + \  \frac{1}{N}\,r^N\,\sin(N\,\theta) \ ,
\end{equation}
here $r=\sqrt{x^2+y^2},\,\theta=\tan^{-1}(y/x)$ are polar coordinates, $p_r$ and $p_{\theta}$ are their canonical conjugate momentum variables, respectively, $m$ is the mass of the particle, $\omega>0$ and $N>0$ is a positive integer number. The phase space is 4-dimensional.

The Hamiltonian (\ref{HN}) possesses a discrete rotational dihedral symmetry $\theta \rightarrow \theta + \frac{2\,\pi}{N}$, and for even $N$ the reflection symmetry $r \rightarrow -r$ is present formally. Furthermore, the generalized Hénon-Heiles systems \eqref{HN} are {\em reversible}. Since ${\cal H}_N$ is a quadratic homogeneous function in the momenta, the Hamilton's equations remain invariant under the time and momenta inversions
\begin{equation}
(r,\,\theta,\, p_r,\,p_\theta,\,t)\ \rightarrow \ (r,\,\theta,\, -p_r,\,-p_\theta,\,-t)\ . 
\end{equation}  

\vspace{0.2cm}

At $N=1,2$, the system is superintegrable and integrable, respectively, whereas for any $N>2$ it exhibits a chaotic motion. In particular, at $N=3$ with $m=\omega=1$ it coincides with the celebrated Hénon-Heiles potential (\ref{HenHei}) at $g=1$.

In Cartesian coordinates, $x=r\,\cos{\theta}$ and $y=r\,\sin{\theta}$, the Hamiltonian (\ref{HN}) reads
\begin{equation}
\begin{aligned}
\label{HNc}
{\cal H}_{{}_N} \ & =  \  \frac{1}{2\,m}\big( p_x^2 \ + \ p_y^2 \big)\ + \ V_N(x,y)
\\ & 
= \frac{1}{2\,m}\big( p_x^2 \, + \, p_y^2 \big)\, + \, 
\frac{m\,\omega^2}{2}\left(x^2+y^2\right) \,  + \, P_N(x,\,y)  \ , 
\end{aligned}
\end{equation}
where again $p_x$ and $p_y$ are the corresponding canonical momentum variables and $P_N(x,y)$ in (\ref{HNc}) is an homogeneous polynomial of order $N$ ($N>1$) in coordinates $x,y$. Explicitly, for the lowest values of $N$ the potential $V_N$ (\ref{VN}) reads
\begin{equation}
\begin{aligned}
\label{pots}
& V_1\ = \ \frac{m\,\omega^2}{2}\left(x^2+y^2\right) \ + \ y 
\\ &
V_2\ = \ \frac{m\,\omega^2}{2}\left(x^2+y^2\right) \ + \ x\,y
\\ &
V_3\ = \ \frac{m\,\omega^2}{2}\left(x^2+y^2\right) \ + \ y\,x^2 \ -\ \frac{1}{3}\,y^3\ 
\\ &
V_4\ = \ \frac{m\,\omega^2}{2}\left(x^2+y^2\right) \ + \ x^3\,y \ -\ y^3\,x\ 
\\ &
V_5\ = \ \frac{m\,\omega^2}{2}\left(x^2+y^2\right) \ + \ y\,x^4 \ -\ 2\,x^2\,y^3\ +\ \frac{1}{5}\,y^5
\\ &
 V_6\ = \ \frac{m\,\omega^2}{2}\left(x^2+y^2\right) \ + \ y\,x^5 \ +\ x\,y^5\ -\ \frac{10}{3}\,y^3\,x^3\ .
\end{aligned}
\end{equation} 

The Hamiltonian (\ref{HNc}) defines a fourth-order autonomous system with variables $(x,y,p_x,p_y)$. It is a dynamical system with polynomial equations of motion.

It is worth noting that under the rescaling $\mathbf{r} \rightarrow \alpha\,\mathbf{r}$ accompanied by the transformations
\begin{equation}
 m \rightarrow \frac{1}{\alpha^{N+2}}\,m \quad , \quad \omega    \rightarrow \alpha^{N}\,\omega \quad , \quad E \rightarrow \alpha^{N}\,E \ ,
\end{equation}
the dynamics of the Hamiltonian (\ref{HN}) is scale-invariant.
\begin{figure}[t!]
	\centering \includegraphics[width=\linewidth]{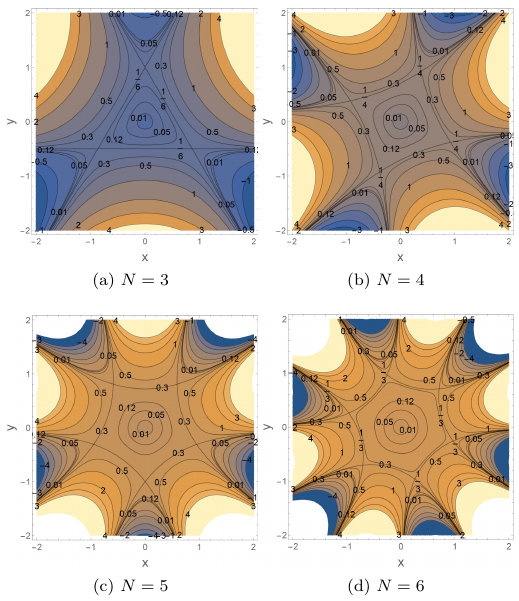}
	\par 
	\caption{Level curves, on the plane $(x,y)$, of the potential $V_N$ (\ref{VN}) as a function of $N$. The corresponding values of the critical energy $E_{\rm crit}=\frac{N-2}{2\,N}$ are displayed as well. The value of the combination $m\,\omega^2=1$ is used in the calculations. }
	\label{LCN3}
\end{figure}

For any $N=1,2,3,\dots$, the time evolution of the  variables $x$ and $y$ takes the form:
\begin{equation}\label{veleq}
\dot x \ = \ \frac{1}{m}\,p_x \quad ; \quad \dot y \ = \ \frac{1}{m}\,p_y \ ,
\end{equation}
thus, they are linear while for the momentum coordinates ($p_x,\,p_y$) the dynamics, $\dot p_x=-\frac{\partial}{\partial\,x}{\cal H}_{{}_N}$ and $\dot p_y=-\frac{\partial}{\partial\,y}{\cal H}_{{}_N}$, is $N-$dependent and non-linear starting from $N=3$. 

In the case $N=1$, making the canonical change of variables $x \rightarrow x$, $y \rightarrow y -\frac{1}{m\,\omega^2}$ with momenta $p_x,p_y$ unchanged, we arrive to the 2D isotropic harmonic oscillator. Thus, the Hamiltonian ${\cal H}_{{}_{N=1}}$, the angular momentum $L_z=x\,p_y-y\,p_x$ as well as $S_{xy}=p_x\,p_y + m^2\,\omega^2\,x\,y$ are three algebraically independent integrals of motion. The system ${\cal H}_{{}_{N=1}}$ and its quantum counterpart are maximally superintegrable \cite{Miller_2013, Escobar-Ruiz_2020}. Similarly, it can be shown that the case $N=2$ corresponds to an integrable system.

\begin{figure*}[t!]
\includegraphics[width=\linewidth]{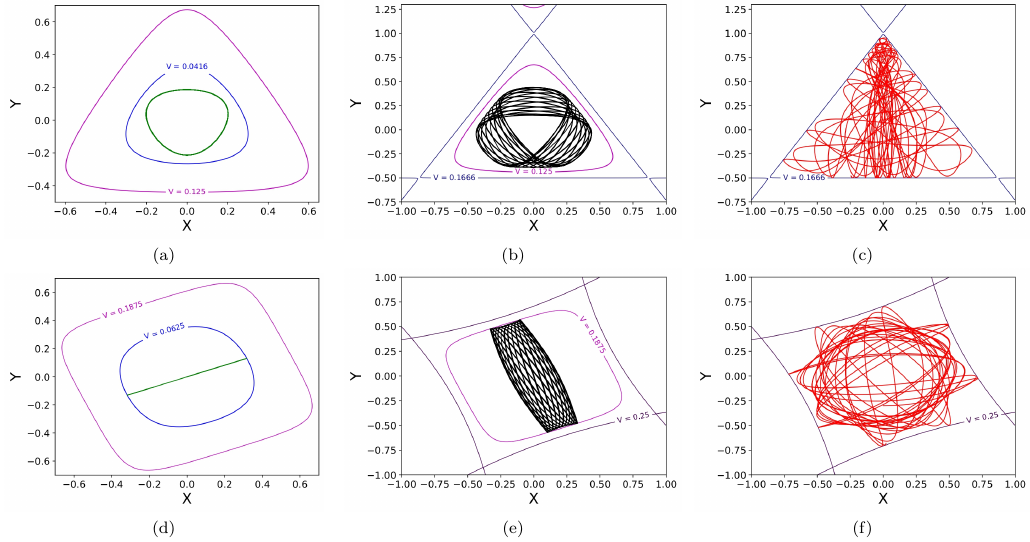}
\caption{Top panel shows level curves on the plane $(x,y)$ for the case $N=3$ and considering the potential $V_3$ (\ref{pots}) with $m\,\omega^2=1$. The corresponding critical energy $E_{\rm crit}=\frac{1}{6}\approx 0.1666$ is indicated as well as certain trajectories: (a) periodic trajectory (green line) with initial conditions $(x_0=0,\,y_0=-0.21341508,\,p_{x_0}=0.17656123,\,p_{y_0}=0)$ and energy $E=\frac{1}{4}E_{c}$, (b) quasi-periodic orbit (dark line) with initial conditions $(x_0=0,\,y_0=0.15,\,p_{x_0}=0.479322,\,p_{y_0}=0)$ and energy $ E=\frac{3}{4}E_{c}$, and (c) a chaotic trajectory (red line) with initial conditions $(x_0=0,\,y_0=0.1,\,p_{x_0}=0.520545,\,p_{y_0}=0.23)$ and energy $E=\frac{99}{100}E_{c}$. Bottom panel shows level curves on the plane $(x,y)$ for the case $N=4$ with potential $V_4$ (\ref{pots}) and parameters $m\,\omega^2=1$. The critical energy $E_{\rm crit}=\frac{1}{4}=0.25$ is indicated as well as certain trajectories: (d) periodic trajectory (green line) with initial conditions $(x_0=0,\,y_0=0,\,p_{x_0}=0.32666587,\,p_{y_0}=0.13523834)$ and energy $E=\frac{1}{4}E_{c}$, (e) quasi-periodic orbit (dark line) with initial conditions $(x_0=0,\,y_0=-0.4,\,p_{x_0}=0.234521,\,p_{y_0}=-0.4)$ and energy $ E=\frac{3}{4}E_{c}$, and (f) a chaotic trajectory (red line) with initial conditions $(x_0=0,\,y_0=-0.3,\,p_{x_0}=0.640312,\,p_{y_0}=0)$ and energy $E=\frac{99}{100}E_{c}$.}
\label{XY_N3}
\end{figure*}

For higher values of $N$, the equations of motion become nonlinear. In particular, at $N=3$, we obtain:
\begin{equation}\label{pdotN3}
\begin{aligned}
&\dot p_x \ = \ -(m\,\omega^2\,x \ + \ 2\,y\,x )
\\ &
\dot p_y \ = \ -(m\,\omega^2\,y \ + \ x^2 \ - \  y^2) \ .
\end{aligned}
\end{equation}

For $N=4$ we have the equations:
\begin{equation}
\begin{aligned}
\label{eqsN4}
&\dot p_x \ = \ -(m\,\omega^2\,x \ + \ 3\,y\,x^2 \ - \ y^3 )
\\ &
\dot p_y \ = \ -(m\,\omega^2\,y \ + \ x^3   \ - \  3\,y^2\,x) \ .
\end{aligned}
\end{equation}

At $N=5$:
\begin{equation}\label{pdotN5}
\begin{aligned}
&\dot p_x \ = \ -( m\,\omega^2\,x \ + \ 4\,y\,x^3 \ - \ 4\,x\,y^3 )
\\ &
\dot p_y \ = \ -( m\,\omega^2\,y \ + \ x^4   \ - \  6\,y^2\,x^2 \ + \ y^4 ) \ .
\end{aligned}
\end{equation}
Hereafter, we put $m=1$ and $\omega=1$. For fixed $N$, there exists a critical value (escape energy)
\begin{equation}
E_{c}(N) \ \equiv \ \frac{N-2}{2\,N}  \ ,
\end{equation}
above which the energy surface ${\cal H}_N=E$ is unbounded, i.e., the trajectories of the system transit from bounded to unbounded motion at $E=E_c$. This critical energy is a monotonously increasing function of $N$. It runs from $-\frac{1}{2}$ at $N=1$ to $\frac{1}{2}$ for $N\rightarrow \infty$. In Figure \ref{LCN3}, for the cases $N=3,4,5,6$ the level curves, on the plane $(x,y)$, of the potential $V_N(x,y)$ (\ref{VN}) are depicted. The dihedral $D_N$ symmetry dictates the regular form of these curves.

\section{Bounded motion}

In this section, for the cases $N=3$ and $N=4$ we explore the bounded dynamics of the Hamiltonian ${\cal H}_N$ (\ref{HNc}) as a function of the energy $E$.

\subsection{Dynamics in configuration space}

In the case of finite motion, the boundaries in configuration space are defined by the equation $V_N = E_c=\frac{N-2}{2\,N}$.
For the values $N=3$ (top panel of Fig. \ref{XY_N3}) and $N=4$ (bottom panel of Fig. \ref{XY_N3}), these boundaries exhibit the 2D dihedral $D_N$ symmetry, as expected. We also present in Fig. \ref{XY_N3} examples of generic periodic, quasi-periodic and chaotic trajectories occurring for $0<E<E_c$. These trajectories were calculated numerically.

\begin{figure*}
    \centering \includegraphics[width=\linewidth]{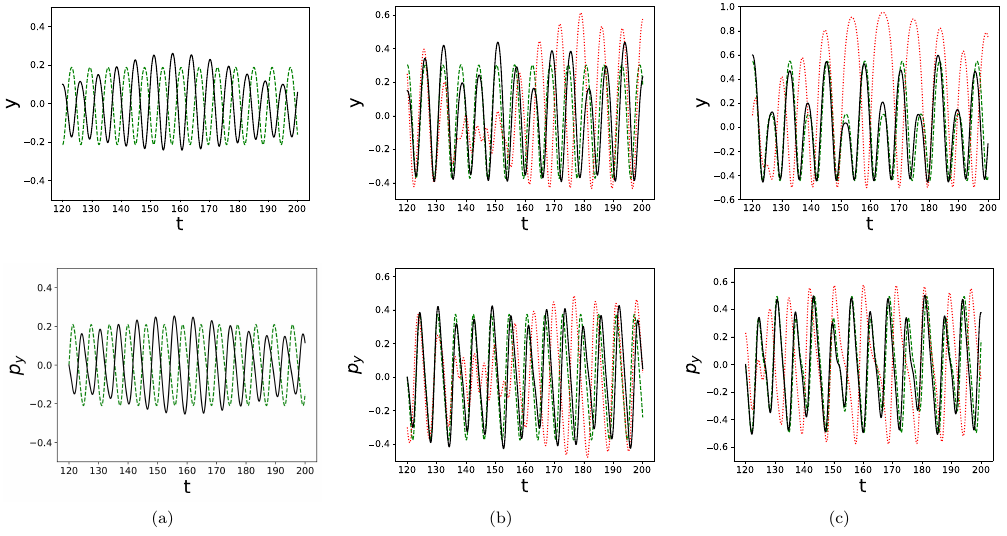}
    \caption{Case $N=3$: time series of the dynamical variables $y(t)$ and $p_{y}(t)$ of ${\cal H}_{{}_{N=3}}$ (\ref{HNc}) for three representative trajectories of the periodic (green dashed), quasi-periodic (black solid) and chaotic (red dotted) motion at different values of the energy. Here $E=\frac{1}{4}E_c$ in (a), $E=\frac{3}{4}E_c$ in (b), whereas $E=\frac{99}{100}E_c$ for (c). The corresponding initial conditions $(x_{0},y_{0},p_{x_{0}},p_{y_{0}})$ are those used for the top panel in Fig. \ref{XY_N3}.}
    \label{N3TS}
\end{figure*}

\begin{figure*}
    \centering
    \includegraphics[width=\linewidth]{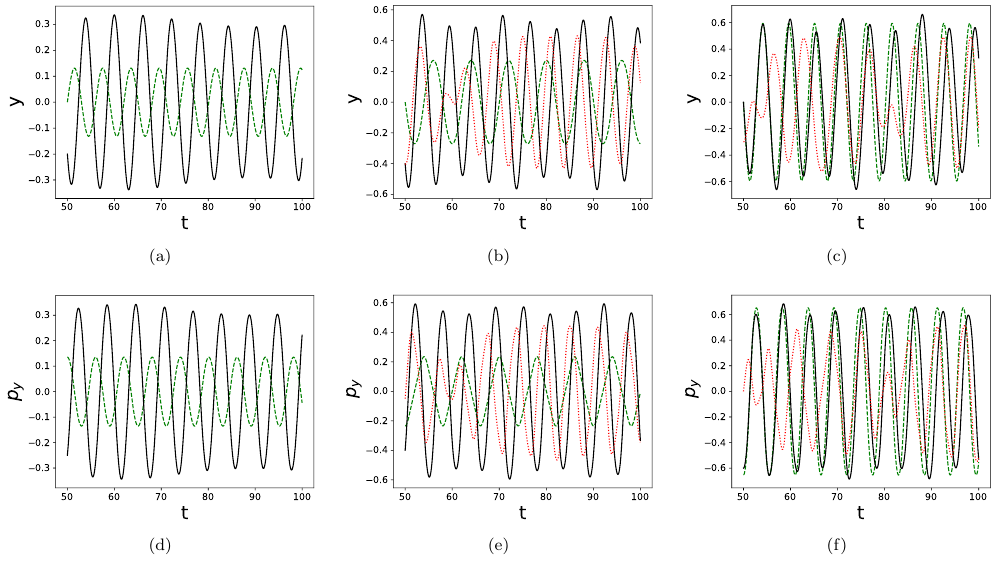}
    \caption{Case $N=4$: time series $y(t)$ and $p_{y}(t)$ of ${\cal H}_{{}_{N=4}}$ (\ref{HNc}) for three representative trajectories of the periodic (green dashed), quasi-periodic (black solid) and chaotic (red dotted) motion at different values of the energy. Here $E=\frac{1}{4}E_c$ in (a), $E=\frac{3}{4}E_c$ in (b), whereas $E=\frac{99}{100}E_c$ for (c). The corresponding initial conditions $(x_{0},y_{0},p_{x_{0}},p_{y_{0}})$ are those used for the bottom panel of Fig. \ref{XY_N3}.}
    \label{N4TS}
\end{figure*}

\begin{figure*}
\centering
\includegraphics[width=\linewidth]{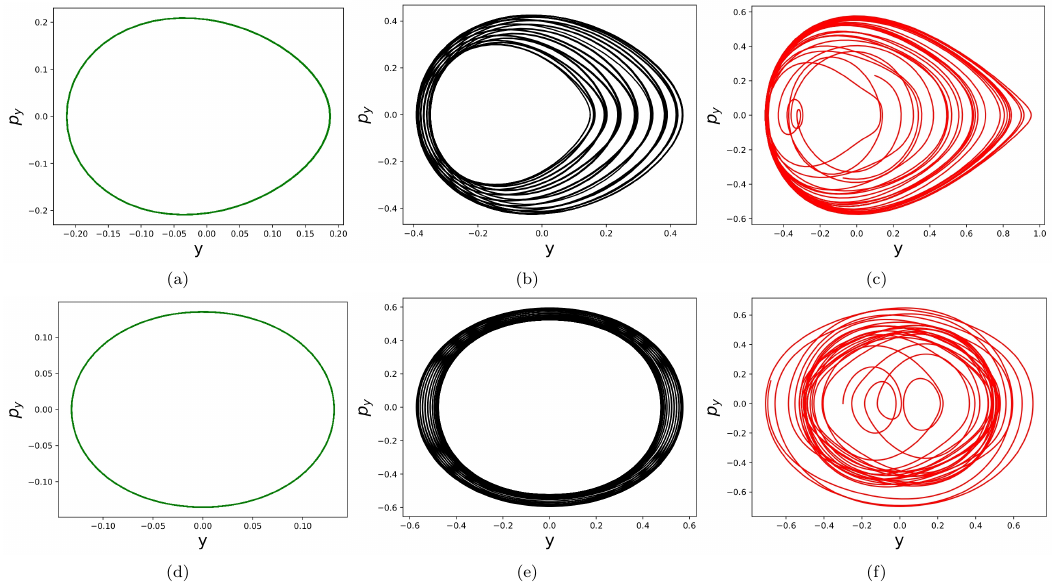}
	\caption{Cases $N=3$ (top) and $N=4$ (bottom): numerical trajectories $p_y=p_y(y)$, corresponding to those displayed in Fig. \ref{XY_N3}, for representative orbits of the periodic (green dashed), quasi-periodic (black solid) and chaotic (red dotted) motion. Here $E=\frac{1}{4}E_c$ in (a) and (d), $E=\frac{3}{4}E_c$ in (b) and (e), while $E=\frac{99}{100}E_c$ for (c) and (f). } 
	\label{N3FD}
\end{figure*}

\begin{figure*}[t!]
\centering
\includegraphics[width=\linewidth]{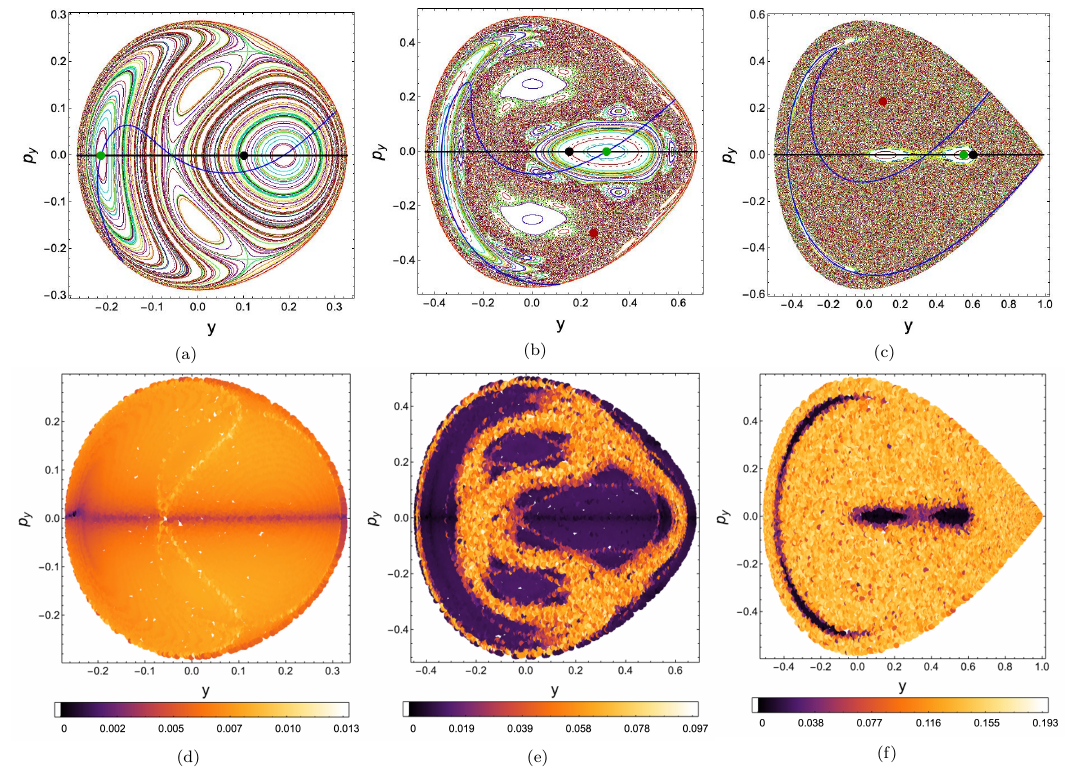}
    \caption{Case $N=3$: (a)-(c) oriented Poincaré sections and symmetry lines $\Gamma_0$ and $\Gamma_2$ (\ref{simlines}), and (d)-(f) the largest Lyapunov exponent on the plane $(y,\,p_y)$, for the Hamiltonian ${\cal H}_{{}_{N=3}}$ (\ref{HNc}) with $m=1$, $\omega=1$ at different values of the energy $E= \frac{1}{4}E_c,\,\frac{3}{4}E_c,\,\frac{99}{100}E_c.$ In this case, the critical energy takes the value $E_c=\frac{1}{6}.$  The (marked) bigger points in green, black and red indicate initial conditions for periodic, quasi-periodic and chaotic trajectories, respectively; they were used to construct Figs. \ref{XY_N3}-\ref{N3FD}}.
    \label{poincareN3}
\end{figure*}

\begin{figure*}[t!]
\centering
\includegraphics[width=\linewidth]{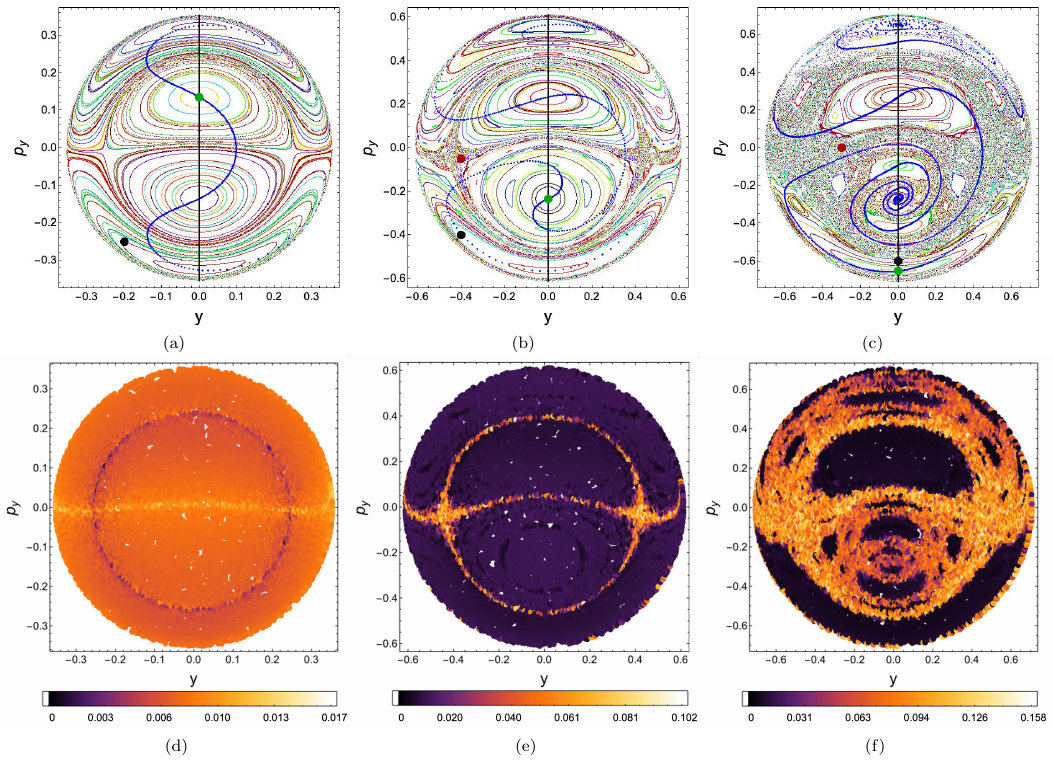}
    \caption{Case $N=4$: (a)-(c) oriented Poincaré sections  and symmetry lines $\Gamma_0$ and $\Gamma_2$ (\ref{simlines}), and (d)-(f) largest Lyapunov exponents on the plane $(y,\,p_y)$, for the Hamiltonian ${\cal H}_{{}_{N=4}}$ (\ref{HNc}) with $m=1$, $\omega=1$ at different values of the energy $E= \frac{1}{4}E_c,\,\frac{3}{4}E_c,\,\frac{99}{100}E_c.$ In this case, the critical energy is $E_c=\frac{1}{4}$, see text. The (marked) bigger points in green, black and red indicate initial conditions for periodic, quasi-periodic and chaotic trajectories, respectively; they were used in Figs. \ref{XY_N3}-\ref{N3FD}}
    \label{poincareN4}
\end{figure*}

\subsection{Time series}\label{subsec:time_series}

Now, for the cases $N=3$ and $N=4$, we display time series (trajectories) obtained by solving numerically the Hamilton's equations of motion. In particular, we plot the data of the dynamical conjugate variables $y(t)$ and $p_y(t)$. For each value of energy $E=\frac{3}{4}E_c,\,\frac{99}{100}E_c$, we present three representative trajectories of the regular (periodic), quasi-periodic and chaotic motion, whereas for $E= \frac{1}{4}E_c$ only periodic and quasi-periodic ones are considered (at this small value of energy there is no significant presence of chaos). The corresponding results for $N=3$ and $N=4$ are shown in Fig. \ref{N3TS} and \ref{N4TS}, respectively. A straightforward Fourier decomposition of the aforementioned time series confirms the nature of these trajectories. For completeness, in Fig. \ref{N3FD} we show the associated dynamical motion in the phase space plane $(y,p_y)$.

%\begin{figure*}[t!]
%\centering  
%\includegraphics[width=\linewidth]{fig6.pdf}
%\par
%\caption{Oriented Poincaré sections, on the plane $(y,\,p_y)$, for the Hamiltonian ${\cal H}_{{}_{N=3}}$ (\ref{HN}) at different values of the energy $E= \frac{1}{4}E_c,\,\frac{3}{4}E_c,\,\frac{99}{100}E_c.$} 
%\label{}
%\end{figure*}

\section{Poincaré sections and Lyapunov exponents}

In this section, we present the oriented Poincaré sections on the $(y,p_y)$ plane, considering $m=1$ and $\omega=1$, for the cases $N=3,4$ as a function of energy $E$. These sections are obtained from the intersection of orbits, associated with specific initial conditions within the state space, with a lower-dimensional subspace ($x=0$, $y$, $p_x$, $p_y$). These sections are transversal planes to the flow of the system. Interestingly, the Poincaré section can be seen as a discretized version of the dynamical system, retaining many characteristics of the original continuous system but operating in a reduced state space \cite{strogatz2018nonlinear, henon1982numerical, tucker2002computing}. 

Complementarily, we analyze the trajectories' nature by calculating the largest Lyapunov exponents. The Lyapunov exponent provides a measure of the average exponential divergence of states that began infinitesimally close. By observing the evolution of nearby trajectories over time, the Lyapunov exponent offers valuable insights into the system's sensitivity to initial conditions, enabling the characterization of chaotic behavior in dynamical systems \cite{sandri1996numerical, young2013mathematical}. The largest Lyapunov exponents can be obtained by computing the double limit: 

\begin{equation}
     \lambda \ = \ \lim_{t \to \infty}\lim_{||\delta \vec{X}(t_0)|| \to 0}\frac{1}{t}\,ln\bigg(\frac{||\delta \vec{X}(t)||}{||\delta \vec{X}(t_0)||}\bigg) \ ,
\end{equation}
where $\delta \vec{X}(t)=\Phi \cdot \delta \vec{X}(t_0)$ with $\vec{X}=[x,y,p_x,p_y]^T$ and $\delta \vec{X}(t_0)$ the perturbed initial state.
$\Phi$ is the fundamental matrix whose elements are calculated by solving $\dot \Phi=J\Phi$ and Hamilton equations, simultaneously. Here, $J$ denotes the Jacobian matrix. 

Figure \ref{poincareN3} and \ref{poincareN4} show oriented Poincaré sections and largest Lyapunov exponents for the Hamiltonians ${\cal H}_{{}_{N=3}}$ and  ${\cal H}_{{}_{N=4}}$, respectively, considering three different energies: $E=\frac{1}{4}E_c$,  $E=\frac{3}{4}E_c$ and $E=\frac{99}{100}E_c$. To draw the Poincaré sections, we develop numerical simulations for $70$ different (random) initial conditions with a simulation time of $7000$ seconds. In contrast, for each value of the energy, Lyapunov exponents plots were constructed by integrating numerically 15000 random initial conditions during 7000 seconds. Computations were performed in Mathematica 12.1. 

%In Figure \ref{poincareN3}(a)-(d), the Poincaré sections present regular dynamics for low energy, $E=\frac{1}{4}E_c$. This is confirmed by Lyapunov exponents extremely close to zero. However, for the cases  $E=\frac{3}{4}E_c$ (Fig. \ref{poincareN3}(b)-(e)) and $E=\frac{99}{100}E_c$ (Fig. \ref{poincareN3}(c)-(f)),  Poincare sections show regions with scattered points whose Lyapunov exponents plots exhibit values different from zero, unveiling the presence of chaos in the system. Importantly, some islands of periodic behaviors are clearly visible in both cases, which are denoted by dark regions. 

In Figure \ref{poincareN3}(a),(d), the Poincaré sections present regular dynamics for low energy, $E=\frac{1}{4}E_c$. This is confirmed by Lyapunov exponents extremely close to zero. However, for the cases  $E=\frac{3}{4}E_c$ (Fig. \ref{poincareN3}(b),(e)),  Poincare sections show regions with scattered points whose Lyapunov exponents plot exhibit values different from zero, unveiling the presence of chaos in the system. Importantly, some islands of periodic behaviors are clearly visible which are denoted by dark regions. Despite the predominance of chaos at higher energy levels, $E=\frac{99}{100}E_c$, these islands of order persist, unveiling the coexistence of regularity and chaos, as displayed in Poincare section and Lyapunov exponents plot of Figs. \ref{poincareN3}(c),(f), respectively. This coalescence of stability allows us to demonstrate the dynamic richness of the system. 

Notably, the Hamiltonian ${\cal H}_{{}{N=4}}$ exhibits regular behavior at low energy, specifically at $E=\frac{1}{4}E_c$, as evidenced by the Poincare sections and Lyapunov exponents plots in Figs. \ref{poincareN4}(a),(d), respectively. When the energy increases up to $E=\frac{3}{4}E_c$, chaotic trajectories begin to be visible although most of the trajectories remain regular as depicted in Figs \ref{poincareN4}(b),(e). Moving to an energy of $E=\frac{99}{100}E_c$,  a majority of initial conditions exhibit chaotic dynamics although some regular regions survive (see Figs \ref{poincareN4}(c),(f)).

\subsection{Symmetry lines }

Additionally, we show the symmetry lines (defined in \cite{PINA1987369}) for the cases of $N=3$ and $N=4$ . These symmetry lines are used for obtaining the initial conditions of periodic orbits.

The generalized Hénon-Heiles systems \eqref{pdotN3}-\eqref{pdotN5} exhibit the property of {\em reversibility}. Specifically, the Hamilton's equations remain unchanged under the inversion of time and momenta, defined as follows:
\begin{equation}\label{ptinv}
 ({\bf q},{\bf p},t)\ \rightarrow \  ({\bf q},-{\bf p},-t)\ .
\end{equation}
Moreover, for even values of $N$, they also possess the symmetry of coordinates and time inversion:
\begin{equation}\label{qtinv}
 ({\bf q},{\bf p},t)\ \rightarrow \ ({\bf -q},{\bf p},-t)\ .
\end{equation}

For these reversible systems, the Poincaré map $T$ can be expressed as the composition of two involutions \cite{PINA1987369}:
\[ T= I_0 T^{-1} I_0 = I_1 I_0\ , \quad  \text{with}\, I_0^2=I_1^2=1 \ , \] 
where $T,\, I_0, I_1:  \Sigma\rightarrow \Sigma$, here $\Sigma$ denotes a Poincaré section. Introducing $I_j=T^j I_0$, the infinite set of transformations $T^j, I_i \, : \Sigma\rightarrow \Sigma$,  with $i,j\in \mathbb{Z}$, forms an infinite discrete group. The fixed points of $I_i$ are referred to as the {\em symmetry lines}:
\begin{equation}
\label{simlines}
\Gamma_i\ = \ \{ {\bf x} \in \Sigma |\, I_i {\bf x}={\bf x}\}, \qquad \text{for all} \,  i\in \mathbb{Z}\ .
\end{equation}

Of particular significance are the fundamental lines $\Gamma_0$ and $\Gamma_1$, they lead to the property
\begin{equation}\label{Gn}
T^i\,\Gamma_j\ = \ \Gamma_{2i+j}\ .  
\end{equation} 
This implies that $\Gamma_0$ generates all symmetry lines with even indices: $\Gamma_{2i}=T^i\Gamma_0$ whilst $\Gamma_1$ produces all symmetry lines with odd indices: $\Gamma_{2i+1}=T^i\Gamma_1$. It can be proven that the intersections of the symmetry lines $\Gamma_i\cap\Gamma_j$ correspond to periodic orbits, with the period being a divisor of the integer $|i-j|$. 

For the Henon-Heiles system \eqref{HN} with $N=3$, the involution $I_0$ is used in the $\Sigma$ plane, along with its fixed points:
\begin{eqnarray}\label{I0ptinv}
I_0(x,y,p_x,p_y,t)&=&(x,y,-p_x,-p_y,-t), \\ \Gamma_0&=&\{ (y,p_y)|p_y=0\}\ .
\end{eqnarray}
For the case of $N=4$, the relevant symmetry is as follows:
\begin{eqnarray}\label{I0qtinv}
I_0(x,y,p_x,p_y,t)&=&(-x,-y,p_x,p_y,-t) \\ \Gamma_0&=&\{ (y,p_y)|y=0\}\ .
\end{eqnarray}

These symmetry transformations are employed to obtain initial conditions for periodic orbits. 

\vspace{0.2cm}

Specifically, for $N=3$ and the symmetry \eqref{I0ptinv}, the periodic orbits correspond to the intersection points of the horizontal line $\Gamma_0$ with $\Gamma_2$ (blue line), Figs. \ref{poincareN3}(a)-(c). Meanwhile, for $N=4$ and the symmetry \eqref{I0qtinv}, the periodic orbits stem from the intersection of the vertical line $\Gamma_0$ with $\Gamma_2$ (blue line), Figs. \ref{poincareN4}(a)-(c).

\section{Reconstruction of the Hamilton's equations from time series}

This section explores the process of reconstructing the Hamilton's equations based only on the time series data computed in Section \ref{subsec:time_series}. First of all, we briefly describe the key elements of the sparse identification of nonlinear dynamical systems (SINDy) algorithm \cite{brunton, dam2017sparse, kaiser2018sparse, lai2019sparse, hoffmann2019reactive}. 

\vspace{0.2cm}

The goal is to estimate quantitatively the accuracy of this data-driven algorithm.  
\subsection{SINDy Algorithm}
\label{Sindyal}
Let $\mathbb{X}\in\mathbb{R}^{n\times4}$ be a matrix in which each column contains $n$ samples of the time series for each of the four dynamical variables $x,\ y,\ p_x $ and $p_y$. Moreover, $\dot{\mathbb{X}}\in\mathbb{R}^{n\times4}$ denotes a matrix whose columns are numerical derivatives of each time series data.  With this in mind, we write the following matrix equation
\begin{equation}\label{eq:ss_data}
    \dot{\mathbb{X}}\ = \ \bm\Theta(\mathbb X)\,\bm \Xi.
\end{equation}
\\
Here, $\bm\Theta(x, y, p_x, p_y)\in\mathbb{R}^{1\times m}$ is the vector of candidate functions that could potentially constitute the reconstructed Hamilton's equations and $\bm\Theta(\mathbb X)\in\mathbb{R}^{n\times m}$ represents the evaluation of the vector of candidate functions for the time series $\mathbb X$. On the other hand, the entries of $\bm\Xi\in\mathbb{R}^{m\times4}$ will be the coefficients of such candidate functions and therefore the unknowns. Since the matrix equation \eqref{eq:ss_data} is overdetermined, we aim to find a solution such that $\bm\Xi$ minimizes $||\dot{\mathbb{X}}\ -\ \bm\Theta(\mathbb X)\,\bm \Xi||$ where $||\cdot||$ denotes the Euclidean norm. 

\vspace{0.1cm}

It is important to note that this least squares solution may not necessarily be sparse, meaning that the majority of its elements are not zero. The work in \cite{brunton} proposed SINDy which is a method to find a sparse solution for the system identification problem, whose results suggest that, with an appropriate choice of candidate functions, a least squares sparse solution will result in a parsimonious model that aligns with the governing equations. Finally, given a sparsity solution obtained via SINDy, the reconstructed model is presented as follows:
\begin{equation}\label{eq:reconstructed}
    \begin{bmatrix}
        \dot x\\
        \dot y\\
        \dot p_x\\
        \dot p_y
    \end{bmatrix}= \bm \Xi^T \bm\Theta(x, y, p_x, p_y)^T.
\end{equation}
Our work uses the Python library PySINDy \cite{de2020pysindy, 2021pysindypackage} to reconstruct the Hamilton's equations from the
time series data computed in Section \ref{subsec:time_series}. For a more comprehensive understanding of the SINDy algorithm, we direct the reader to \cite{brunton}.

\subsection{Results and Discussion}

As previously mentioned in subsection \ref{Sindyal}, concerning the reconstruction of state equations from the time series, it is necessary to choose candidate functions that will comprise the reconstructed model. A physical and mathematically based selection is fundamental to obtain optimal results. Since in this case we know the exact model (\ref{HNc}), for each state equation $\dot{x},\,\dot{y},\,\dot{p}_x,\,\dot{p}_y$, a suitable option is to take a polynomial function of degree $K>N-1$ in the variables $x,\,y,\,p_x,\,p_y$ whose $C$-coefficients will be the decision variables (the unknowns). For example, the most general polynomial ansatz of order $K$ is of the form

\begin{equation}
\label{xans}
   \dot x  \ = \ \sum_{0 \leq n_1+n_2+n_3+n_4=K} C_{\dot x}(n_1,n_2,n_3,n_4)\,x^{n_1}\,\,y^{n_2}\,p_x^{n_3}\,p_y^{n_4}\ ,
\end{equation}
being the $C$-constants real parameters. Similar expressions for the remaining three state equations are employed.

\vspace{0.2cm}

In order to determine the accuracy of SINDy, let us introduce the error parameter $\Delta C_\tau(\sigma)$ as follows
\begin{equation}
\label{DC}
    \Delta C_\tau(\sigma) \ = \  \frac{C_\tau^{\rm exact}(\sigma) \ - \ C_\tau^{\rm SINDy}(\sigma)}{C_\tau^{\rm exact}(\sigma)}\ ,
\end{equation}
here, $C_\tau^{\rm exact}(\sigma)$ represents the coefficient of the monomial $\sigma$ in the state equation $\tau$. 

\vspace{0.2cm}

For the sake of clarifying our notation, let us examine the scenario where $N = 3$. The exact state equation for the derivative of $p_x$ with respect to time is 
$
\dot p_x \ = \ -x \ - \ 2\,y\,x 
$.
Thus, $\tau=\dot p_x$, $\sigma=x, x\,y$, and the non-zero coefficients are $C_{\dot p_x}^{\rm exact}(x)=-1$ and $C_{\dot p_x}^{\rm exact}(x\,y)=-2$ only. Similarly, $C_\tau^{\rm SINDy}(\sigma)$ will denote the coefficient of the monomial $\sigma$ in the (approximate) state equation for $\tau$ constructed by the SINDy algorithm. 

\vspace{0.1cm}

Note that the constant (\ref{DC}), which measures the precision of SINDy, depends on the value $dt$ of the subinterval of time used in the numerical integration of the equations of motion. In general, for sufficiently small values of $dt$ and taking solely the set of data corresponding to a single trajectory, the algorithm of SINDy reproduces exactly the same functional form of the original Hamilton's equations. Naturally, as we increase the value of $dt$, extra terms in the reconstructed equations occur. We define a critical time step $dt_c$ as the smallest value of $dt$ for which an extra term appears in any of the state equations. 

\vspace{0.1cm}

Summarizing, $N$ is fixed and for each value of energy $E=\frac{1}{4}E_c,\,\frac{3}{4}E_c,\,\frac{99}{100}E_c$, we select a trajectory to generate the data (a numerical time series) that (together with polynomial candidate functions) will be used as an input to reconstruct via SINDy the (state) Hamilton's equations. Afterwards, we calculate the corresponding value of $\Delta C_\tau(\sigma)$ as a function of the step $dt$. This process is repeated in the same manner for a periodic, quasi-periodic and chaotic trajectory. It will allow us to study the changes in accuracy when the motion transits from a regular (periodic) to a chaotic behavior. 

\vspace{0.1cm}

\begin{figure}[H]
\centering
\includegraphics[width=0.9\linewidth]{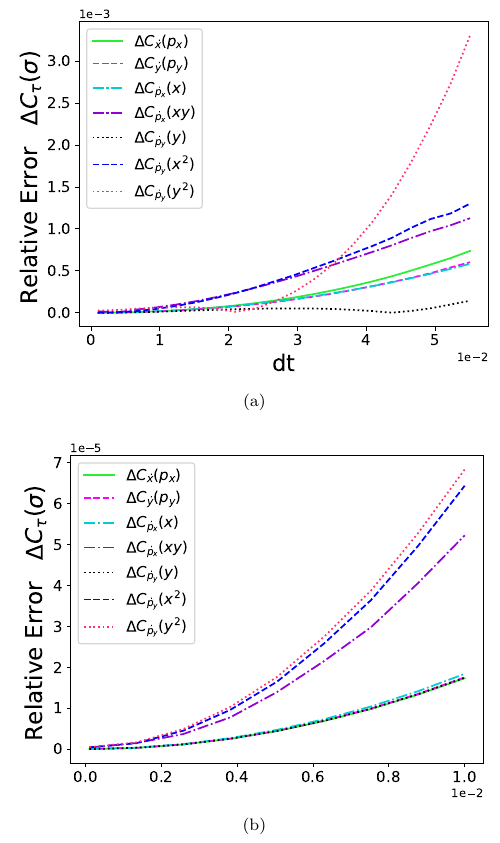}
	\par 
	\caption{Case $N=3$. The relative error $\Delta C_\tau(\sigma)$, see (\ref{DC}), of the coefficients obtained by SINDy for the reconstructed Hamilton's equations; (a) periodic trajectory 
and (b) quasi-periodic case. The critical time integration increments (see text) are $dt_{c}=0.01763$ and $dt_{c}=0.02628$, respectively. These trajectories correspond to the energy value $E=\frac{1}{4}E_{c}$.} 
\label{CN3A}
\end{figure}

\begin{figure}[H]
\centering
\includegraphics[width=0.9\linewidth]{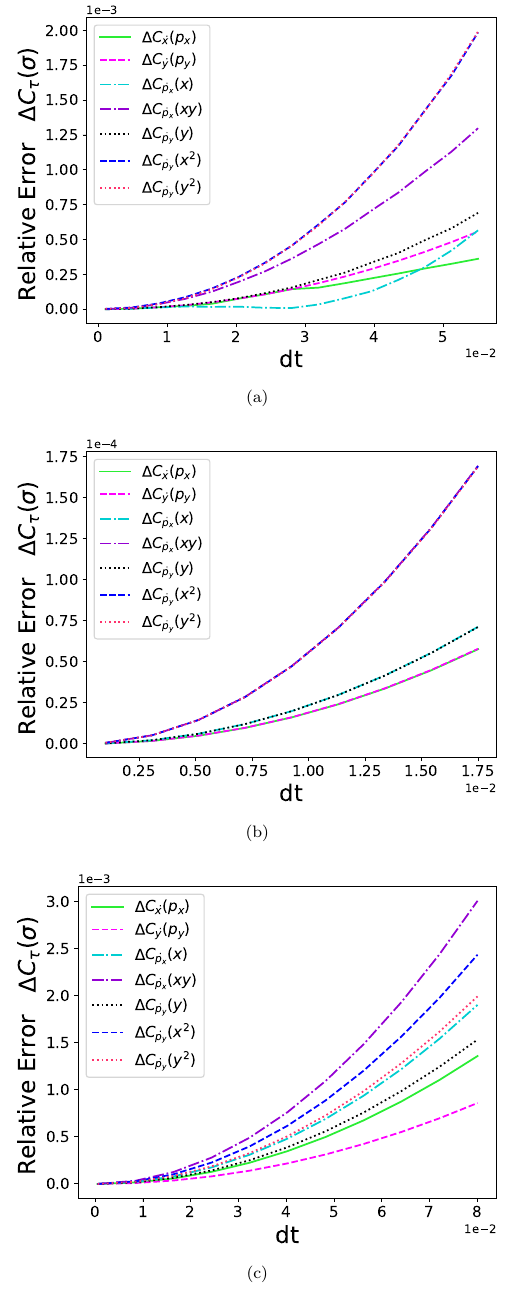}
	\caption{Case $N=3$. The relative error $\Delta C_\tau(\sigma)$, see (\ref{DC}), of the coefficients obtained by SINDy for the reconstructed Hamilton's equations; (a) periodic trajectory, (b) quasi-periodic and (c) chaotic motion. The critical time increments are $dt_{c}=0.02120$, $dt_{c}=0.02406$ and $dt_{c}=0.08443$, respectively. These trajectories correspond to the same energy value $E=\frac{99}{100}E_{c}$.} 
\label{CN3B}
\end{figure}

\begin{figure}[H]
\centering
\includegraphics[width=0.9\linewidth]{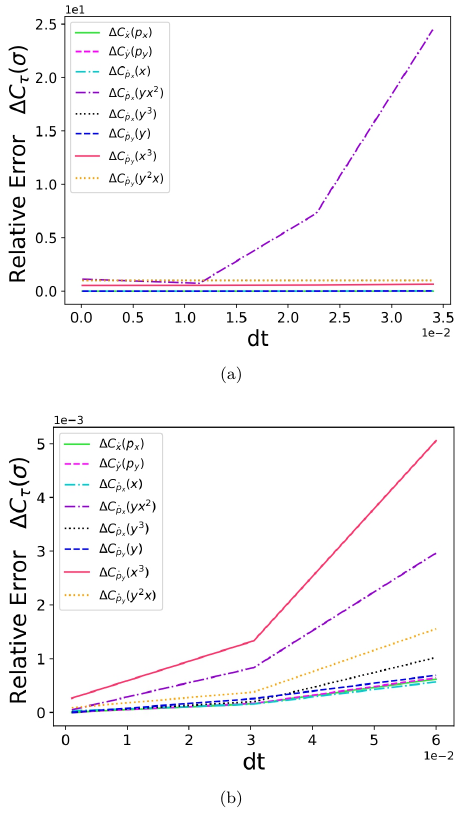}
	\par 
	\caption{Case $N=4$. The relative error $\Delta C_\tau(\sigma)$, see (\ref{DC}), of the coefficients obtained by SINDy for the reconstructed Hamilton's equations; (a) periodic trajectory 
and (b) quasi-periodic case. These trajectories correspond to the energy value $E=\frac{1}{4}E_{c}$. The critical time increment for the periodic trajectory is zero, but for the quasi-periodic case is given by $dt_{c} = 0.0305$. } 
	\label{CN4A} 
\end{figure}

\begin{figure}[H]
\centering
\includegraphics[width=0.9\linewidth]{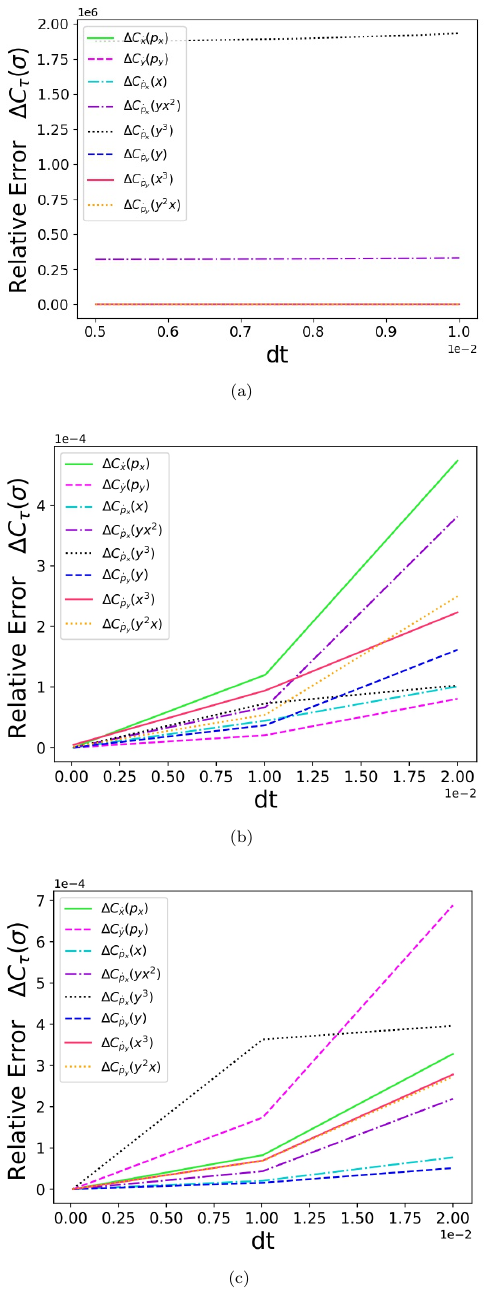}
    \caption{Case $N=4$. The relative error $\Delta C_\tau(\sigma)$, see (\ref{DC}), of the coefficients obtained by SINDy for the reconstructed Hamilton's equations; (a) periodic trajectory, (b) quasi-periodic and (c) chaotic motion. For (a) the critical time increment is zero whereas for (b) and (c), they are $dt_{c}=0.0136$ and $dt_{c}=0.01005$,  respectively. These three trajectories correspond to the energy value $E=\frac{3}{4}E_{c}$.}
    \label{CN4B}
\end{figure}

Concrete results, taking $K=N$ and the threshold $\Gamma=0.05$, for the cases $N=3$ and $N=4$ are presented in Figs. \ref{CN3A}, \ref{CN3B} and Figs. \ref{CN4A}, \ref{CN4B}, respectively. In general, most of the curves $\Delta C_\tau(\sigma)$ \textit{vs} $dt$ are smooth increasing positive functions vanishing as $dt \rightarrow 0$. This implies that, regardless of the nature of motion (periodic or chaotic), the SINDy algorithm recovers the correct underlying dynamical equations of the system where a single trajectory was employed in the reconstruction. 

\vspace{0.1cm}

Essentially, at fixed ($N,E,dt$), the parameter $\Delta C_\tau(\sigma)$ takes its minimum value for a periodic trajectory whereas for a chaotic one it reaches its maximum. Nevertheless, if $dt$ is sufficiently small $\approx10^{-4}$, the presence of chaos does not influence the accuracy provided by SINDy.

\subsubsection{On approximate analytical periodic solutions}

More importantly, as a by product of the present analysis we observed an interesting feature of SINDy that eventually leads to accurate approximate analytical expressions for the periodic trajectories. Therefore, even in the case when the exact model is known and, thus, no reconstruction of the dynamical equation is really needed, SINDy can still be exploited to derive approximate solutions explicitly, see below.   

\vspace{0.1cm}

Taking $N=4$, the energy $E=\frac{1}{4}E_{c}$, the small integration step $dt=2\times 10^{-5}$, and the initial conditions used in Fig. \ref{XY_N3} (d), see also the bigger green point in Fig. \ref{poincareN4} (a), which correspond to a periodic trajectory, the reconstructed Hamilton's equations obtained with SINDy are given by: 
\begin{equation}
\begin{aligned}
\label{SINDyN4}
& \dot x_{_{\scriptsize\text{SINDy}}}\ = \ p_{x}
\\&
\dot y_{_{\scriptsize\text{SINDy}}}\ = \ p_{y}
\\&
\dot p_{x_{\scriptsize\text{SINDy}}}   = \  0.00322\,y \ - \ 1.00132\,x \ - \ 0.75788\,x\,y^2
\\ &  \hspace{1.3cm} 
+ \ 2.81860\,y\,p_x\,p_y \ - \ 1.88935\,y\,x^2 \ - \ 1.19774\,y\,p_x^2 
\\ & \hspace{1.3cm}  
  - \ 1.16865\,x\,p_x\,p_y \ - \ 0.25909\,x^3
 \   + \ 0.49646\,x\,p_x^2
\\ &
\dot p_{y_{\scriptsize\text{SINDy}}}   =  -\,y  \ + \ 2.48528 \,y\,x^2 \  - \ 0.00042\,x\,p_x\,p_y
\\ & \hspace{1.3cm} 
 - \ 1.51478\,x^3
\ + \ 0.00012\,x\,p_x^2
\ ,
\end{aligned}
\end{equation}
here the coefficients were computed with an accuracy of $\approx 3-4$ significant digits. Since we know the exact model, the rhs of the above equations can be decomposed into the sum of the exact part plus extra terms, namely
\begin{equation}
\label{sindecom}
    S_{\tau,\,\tiny\text{SINDy}} \ = S_{\tau,\,{\rm exact}} \ + \ S_{\tau,\,{\rm extra}} \ .
\end{equation}
Clearly, no extra terms occur for $\dot x_{_{\scriptsize\text{SINDy}}}$ and $\dot y_{_{\scriptsize\text{SINDy}}}$. However, in the case of $\dot p_{x_{\scriptsize\text{SINDy}}}$ and $\dot p_{y_{\scriptsize\text{SINDy}}}$ they read
\begin{equation}
\begin{aligned}
\label{eqsN4}
& 
S_{\dot{p}_{x},\text{extra}} \ =  0.00322\,y \ - \ 0.00132\,x \ - \ 0.75788\,x\,y^2 \ - \ y^3
\\ &  \hspace{1.3cm} 
+ \ 2.81860\,y\,p_x\,p_y \ + \ 1.11065\,y\,x^2 \ - \ 1.19774\,y\,p_x^2 
\\ & \hspace{1.3cm}  
 - \ 1.16865\,x\,p_x\,p_y \ - \ 0.25909\,x^3
 \  + \ 0.49646\,x\,p_x^2
\\ &
S_{\dot{p}_{y},\text{extra}}\ = \  2.48528\,y\,x^2 \ - \ 3\,x\,y^2 
  - \ 0.00042\,x\,p_x\,p_y 
\\ & \hspace{1.3cm}
  - \ 0.51478\,x^3
\ + \ 0.00012\,x\,p_x^2
\ ,
\end{aligned}
\end{equation}
respectively.

\begin{figure}[H]
\centering
\includegraphics[width=\linewidth]{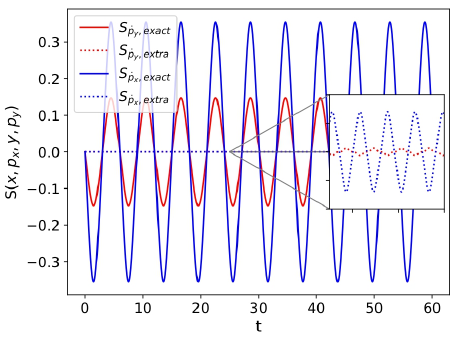}
	\par 
	\caption{Case $N=4$: the time series of the two terms in (\ref{sindecom}), the exact and extra parts, obtained by SINDy. These terms were evaluated over the periodic trajectory displayed in Fig. \ref{XY_N3} (d). The variation of the extra terms is of order $\leq 10^{-7}$.} 
	\label{NeglibS}
\end{figure}

%\bigskip

Now, evaluated over the original numerical trajectory with initial conditions of Fig. \ref{XY_N3} (d), the contribution of $S_{\dot{p}_{x},\text{extra}}$ and $ S_{\dot{p}_{y},\text{extra}}$ above is effectively zero!, see Fig. \ref{NeglibS}. 
In particular, from the second equation in (\ref{eqsN4}) neglecting the terms with coefficients of order $\approx 10^{-4}$, we arrive to the approximate analytical relation 
\begin{equation}
    y \ = \ 0.41\,x\ , \nonumber
\end{equation}
which corresponds with the periodic trajectory in configuration space displayed in Fig. \ref{XY_N3} (d). Hence, for arbitrary small values of $dt$, SINDy can generate spurious terms in the reconstruction of the governing equations of motion.

\subsubsection{$N=3$}

\begin{figure}[H]
\centering
\includegraphics[width=0.9\linewidth]{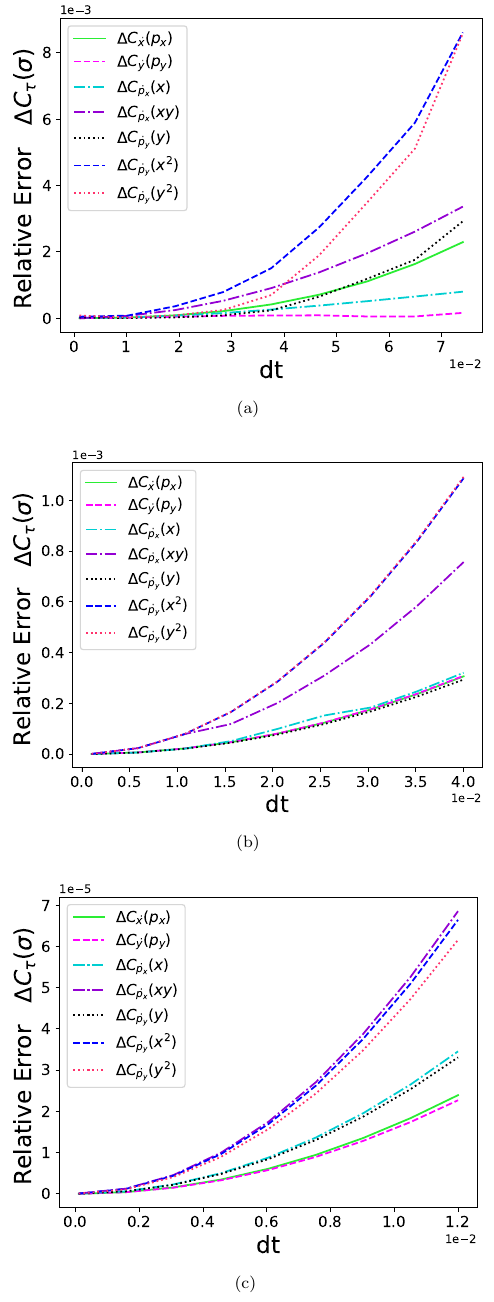}
    \caption{Case $N=3$. The relative error $\Delta C_\tau(\sigma)$, see (\ref{DC}), of the coefficients obtained by SINDy for the reconstructed Hamilton's equations; (a) periodic trajectory, (b) quasi-periodic and (c) chaotic motion. The critical time increments are $dt_{c}=0.01925$, $dt_{c}=0.02726$ and $dt_{c}=0.13303$, respectively. These trajectories correspond to the same energy value $E=\frac{3}{4}E_{c}$.}
\end{figure}

For instance, the periodic trajectory $y=y(x)$ occurring at $N=3$ in Fig. \ref{XY_N3} (d)  does not admits an accurate representation in the polynomial basis (candidate functions) employed by SINDy. Thus, in that case the reconstruction of the dynamical equations is indeed correct (no extra terms appear).  

\vspace{0.1cm}

\textbf{Remark:} \textit{For a given Hamiltonian system, the SINDy algorithm together with an educated election of the candidate functions (a basis) can be used to find approximate analytical expressions of periodic trajectories.} 

\subsubsection{$N=4$}

\section{Conclusions}

In this study, we defined a Hénon-Heiles like potential $V_N=\frac{1}{2}\,m\,\omega^2\,r^2 +   \frac{1}{N}\,r^N\,\sin(N\,\theta)$, $N \in {\mathbb Z}^+$. For arbitrary $N>2$, in Cartesian coordinates $(x,y)$, it corresponds to a polynomial function of degree $N$ possessing a $D_N$ dihedral symmetry. For $E<E_c$, the motion is finite and the system displays a chaotic behavior. A systematic characterization of the dynamics of the system as a function of $N$ and the energy $E$ was carried out. In particular, the explicit calculation of symmetry lines allowed us to compute periodic trajectories with high accuracy. In the second part of the paper, this model was employed to test SINDy which is a method to find the underlying governing (Hamilton's) equations using only the time series data of the trajectories obtained numerically. In general, for sufficiently small time-step of the data, SINDy provides convergent results (even for chaotic trajectories) to the exact functional form of the dynamical equations. Otherwise, the chaotic behaviour decreases the accuracy of the method. Interestingly, it was pointed out that when the corresponding series time data admits an accurate finite expansion in the corresponding basis of candidate functions, then this data-driven algorithm cannot recover the equations of motion correctly. Notably, for a given Hamiltonian, that is to say, when a reconstruction of the governing equations is not needed, SINDy can also be exploited to obtain approximate analytical expressions for the periodic trajectories which is, in general, a non-trivial task. For future work, among the interesting open questions we can mention the quantum counterpart of the model as well as the extension to the case with rational $N$. We plan to study the use of the SINDy algorithm and other methods alike to find symmetries (integrals of motion) of generic Hamiltonian systems as well.

\subsection*{Credit authorship contribution statement}

\textbf{A. M. Escobar-Ruiz} Conceptualization (lead); Formal analysis (equal); Writing original draft
(lead); Writing review/editing (equal), \textbf{L. Jiménez-Lara} Conceptualization (equal); Formal analysis (lead); Investigation (equal);  Writing original draft (equal); Writing review/editing (equal), \textbf{P. M. Juárez-Florez} Formal analysis (equal); Investigation (equal); Validation (supporting), \textbf{F.
Montoya-Molina} Formal analysis (equal); Investigation (equal); Validation (lead), \textbf{J. Moreno-Sáenz} Conceptualization (equal);  Formal analysis (equal); Writing original draft
(equal); Writing review/editing (equal), \textbf{M. A. Quiroz-Juarez} Conceptualization (equal); Formal analysis (equal); Investigation (lead); Writing original draft
(equal); Writing review/editing (equal); Validation (equal).

\subsection*{Declaration of competing interest}

The authors declare that they have no known competing financial interests or personal relationships that could have appeared to influence the work reported in this paper.

\subsection*{Data availability}

No data was used for the research described in the article.

\section*{Acknowledgments}
This work was supported by Consejo Nacional de Humanidades, Ciencias y Tecnologías (CONAHCyT) of Mexico under Grant CF-2023-I-1496. AMER is thankful to CFATA (UNAM) for kind hospitality extended to him where this work was mostly completed. MAQJ would  like  to  thank  the  support  from Dirección General de Asuntos del Personal Académico (DGAPA)-National Autonomous University of Mexico under Project UNAM-PAPIIT TA101023.

%\bigskip

%\newpage
%\nocite{*}
\bibliographystyle{ieeetr}
\bibliography{Manuscript}% Produces the bibliography via BibTeX.Manuscript

\end{document}